\documentclass[11pt,english]{amsart}
\usepackage[T1]{fontenc}
\usepackage{geometry}
\usepackage{amssymb}

\makeatletter
 \theoremstyle{plain}    
 \newtheorem{thm}{Theorem}[section]
 \numberwithin{equation}{section} 
 \numberwithin{figure}{section} 
 \theoremstyle{plain}
 \theoremstyle{definition}
 \newtheorem{defn}[thm]{Definition}
 \theoremstyle{remark}
 \newtheorem{rem}[thm]{Remark}
 \theoremstyle{plain}    
 \newtheorem{prop}[thm]{Proposition} 
 \theoremstyle{plain}    
 \newtheorem{cor}[thm]{Corollary} 

\newcommand\rest{\hbox{\raise.17em\hbox{$ |\kern-.2em$}\lower.23em\hbox{$-$}}}
\newcommand\avint{\hbox{\hbox{$\displaystyle \int$}\hbox{\kern-.9em{$-$}}}}
\newcommand\smavint{\hbox{\hbox{$\int$}\hbox{\kern-.75em{$-$}}}}

\newcommand{\R}{\mathbb{R}}

\newcommand{\Z}{\mathbb{Z}}

\usepackage{babel}
\makeatother
\begin{document}

\title{Regularity of Volume-Minimizing Flows on 3-Manifolds}

\author{David L. Johnson and Penelope Smith}

\date{\today}

\address{Department of Mathematics\\
Lehigh University\\
Bethlehem, Pennsylvania 18015-3174}

\email{dlj0@lehigh.edu, ps02@lehigh.edu}

\begin{abstract}
In \cite{js-1,js-partial-regularity,dim-3} the authors characterized
the singular set (discontinuities of the graph) of a volume-minimizing
rectifiable section of a fiber bundle, showing that, except under
certain circumstances, there exists a volume-minimizing rectifiable
section with the singular set lying over a codimension-3 set in the
base space. In particular, it was shown that for 2-sphere bundles
over 3-manifolds, a minimizer exists with a discrete set of singular
points. 

In this article, we show by analysis of the characterizing horizontal
tangent cone, or $h$-cone, that for a 2-sphere bundle over a compact
3-manifold, such a singular point cannot exist. As a corollary, for
any compact 3-manifold, there is a $C^{1}$ volume-minimizing one-dimensional
foliation. In addition, this same $h$-cone analysis is used to show
that the examples, due to Sharon Pedersen \cite{Pedersen}, of potentially
volume-minimizing rectifiable sections (rectifiable foliations) of
the unit tangent bundle to $S^{2n+1}$ are \textit{not}, in fact,
volume minimizing.
\end{abstract}



\maketitle

\section{Introduction}

In \cite{GZ}, Herman Gluck and Wolfgang Ziller asked which one-dimensional,
transversely oriented  foliation $\mathcal{F}$ (called a \textit{flow})
on an odd-dimensional round  sphere is best-organized, in the sense
that the image of the  natural section $\xi:M\rightarrow T_{1}(M)$
of the unit tangent bundle,  whose value at $x$ is the unit tangent
vector of the leaf of $\mathcal{F}$ through $x$ consistent with
the orientation of $\mathcal{F}$, has smallest $n$-dimensional Hausdorff
measure. 

Their work was in part an effort to interpret the behavior of the
Hopf fibration of the three-sphere, and indeed they were able to show
that the Hopf fibration did minimize the volume. Specifically, they
were able to show that there is a three-form on $T_{1}(S^{3})$ which
calibrates the fibers of the Hopf fibrations on $S^{3}$, thus those
foliations have the least volume of all such flows on the round three-sphere.
However, in higher dimensions the Hopf fibrations are not volume-minimizing,
and it is likely that volume-minimizing flows on these manifolds are
singular. In her thesis \cite{Pedersen}, Sharon Pedersen illustrated
a stable, singular foliation which has much less mass than the Hopf
fibration of $S^{5}$. 

The purpose of the present work is to show that the regularity of
Gluck and Ziller's volume-minimizing flow on $S^{3}$ is a special
case of a theorem that there is a regular ($C^{1}$ as a foliation)
volume-minimizing flow on compact, oriented 3-manifold. Similarly,
there are volume-minimizing sections of the unit tangent bundle (or
other $(n-1)$--sphere bundles over $n$-dimensional manifolds) without
isolated poles. As a corollary result, it will follow that Pedersen's
currents are \textit{not} volume-minimizing among rectifiable sections
of $T_{1}(S^{2n+1})$.

\subsection{Volume of Foliations.}

The \textit{volume} of a one-dimensional foliation $\mathcal{F}$
on a compact manifold $M$ can be computed in terms of the Gauss map
$\xi:M\rightarrow T_{1}(M)$ defined by mapping $x$ to a unit vector
$\xi(x)$ tangent to $\mathcal{F}$ at $x$, which can be chosen consistently
if $\mathcal{F}$ is oriented. The formula is given as:\begin{eqnarray*}
\mathcal{V}(\xi) & = & \int_{M}\sqrt{1+\left\Vert \nabla\xi\right\Vert ^{2}+\dots+\left\Vert \nabla\xi^{\wedge(n-1)}\right\Vert ^{2}}\, dV_{M}\end{eqnarray*}

where the vector wedge is interpreted by\[
\nabla\alpha\wedge\nabla\beta(X,Y):=\frac{1}{2}(\nabla_{X}\alpha\wedge\nabla_{Y}\beta-\nabla_{Y}\alpha\wedge\nabla_{X}\beta),\]
etc., so that \[
(\nabla\xi)^{\wedge k}(X_{1},\dots,X_{k})=\nabla_{X_{1}}\xi\wedge\dots\wedge\nabla_{X_{k}}\xi.\]
The sum is taken over wedges of order up to $n-1$ since the fiber
($S^{n-1}$) is $(n-1)$-dimensional.  Although this is precisely
the $n$-dimensional Hausdorff measure of the image, which is the
mass of the rectifiable current representing the Gauss map as a current
in $T_{1}(M)$, this description has certain advantages.

This definition can be extended to sections $\sigma$ of any smooth
fiber-bundle $B\rightarrow M$ with compact fiber $F$, as defined
in \cite{js-1,js-partial-regularity}. The volume functional is essentially
the same, except that the highest-degree term in the square root is
the minimum of the dimension of $M$ or that of the fiber, \begin{eqnarray*}
\mathcal{V}(\sigma) & = & \int_{M}\sqrt{1+\left\Vert \nabla\sigma\right\Vert ^{2}+\dots+\left\Vert \nabla\sigma^{\wedge n}\right\Vert ^{2}}\, dV_{M},\end{eqnarray*}
with terms $\left\Vert \nabla\sigma^{i}\right\Vert ^{2}$ being 0
for $i>dim(F)$. The results of this article will apply equally to
any $S^{n-1}$-bundle over a compact, oriented $n$-manifold $M$,
but the main impetus of the research came out of the original question
regarding foliations.

\section{\label{sec:Rectifiable-Sections.}Rectifiable Sections.}

Let $B$ be a Riemannian fiber bundle with compact fiber $F$ over
a Riemannian $n$-manifold $M$, with projection $\pi:B\rightarrow M$
a Riemannian submersion. $F$ is a $j$-dimensional compact Riemannian
manifold. Following \cite{MS}, $B$ embeds isometrically in a vector
bundle $\pi:E\rightarrow M$ of some rank $k\geq j$, which has a
smooth inner product $<\:,\:>$ on the fibers, compatible with the
Riemannian metric on $F$. The inner product defines a collection
of connections, called \textit{metric connections}, which are compatible
with the metric. Let a metric connection $\nabla$ be chosen. The
connection $\nabla$ defines a Riemannian metric on the total space
$E$ so that the projection $\pi:E\rightarrow M$ is a Riemannian
submersion and so that the fibers are totally geodesic and isometric
with the inner product space $E_{x}\cong\R^{k}$ \cite{Sasaki}, \cite{j}. 

We will be using multiindices $\alpha=(\alpha_{1},\dots,\alpha_{n-l})$,
$\alpha_{i}\in\{1,\dots,n\}$ with $\alpha_{1}<\cdots<\alpha_{n-l}$,
over the local base variables, and $\beta=(\beta_{1},\dots,\beta_{l})$,
$\beta_{j}\in\left\{ 1,\dots,k\right\} $ with $\beta_{1}<\cdots<\beta_{l}$,
over the local fiber variables (we will at times need to consider
the vector bundle fiber, as well as the compact fiber $F$; which
is considered will be clear by context). The range of pairs $(\alpha,\beta)$
is over all pairs satisfying $|\beta|+|\alpha|=n$, where $\left|\left(\alpha_{1},\dots,\alpha_{m}\right)\right|:=m$.
As a notational convenience, denote by $n$ the $n$-tuple $n:=(1,\ldots,n)$,
and denote the null $0$-tuple by $0$.

\begin{defn}
An $n$-dimensional current $T$ on a Riemannian fiber bundle $B$
over a Riemannian $n$-manifold $M$ locally, over a coordinate neighborhood
$\Omega$ on $M$, decomposes into a collection, called \emph{components},
or \emph{component currents} \emph{of} $T$, with respect to the bundle
structure. Given local coordinates $(x,y)$ on $\pi^{-1}(\Omega)=\Omega\times\R^{k}$
and a smooth $n$-form $\omega\in E^{n}(\Omega\times\R^{k})$, $\omega:=\omega_{\alpha\beta}dx^{\alpha}\wedge dy^{\beta}$,
define auxiliary currents $E_{\alpha\beta}$ by $E_{\alpha\beta}(\omega):=\int\omega_{\alpha\beta}d\left\Vert T\right\Vert $,
where $\left\Vert T\right\Vert $ is the measure $\theta\mathcal{H}^{n}\rest Supp(T)$,
with $\mathcal{H}^{n}$ Hausdorff $n$-dimensional measure in $\Omega\times\R^{k}$
and $\theta$ the multiplicity of $T$ \cite[pp 45-46]{Morgan}. The
\emph{component currents} of $T$ are defined in terms of \emph{component
functions} $t_{\alpha\beta}:\Omega\times\R^{k}\rightarrow\R$ and
the auxiliary currents, by:\[
\left.T\right|_{\pi^{-1}(\Omega)}:=\left\{ T_{\alpha\beta}\right\} :=\left\{ t_{\alpha\beta}E_{\alpha\beta}\right\} ,\]
where $t_{\alpha\beta}$ is defined below. The component functions
$t_{\alpha\beta}:\pi^{-1}(\Omega)\rightarrow\R$ completely determine
the current $T$, and the pairing between $T$ and an $n$-form $\omega\in E^{n}(E)\rest\Omega\times\R^{k}$
is given by:\[
T(\omega):=\int_{\Omega\times\R^{k}}\sum_{\alpha\beta}t_{\alpha\beta}\omega_{\alpha\beta}d\left\Vert T\right\Vert .\]

\end{defn}

\begin{defn}
\label{def:quasi-section}A bounded current $T$ in $E$ is a \emph{(bounded)
quasi-section} if, for each coordinate neighborhood $\Omega\subset M$, 
\begin{enumerate}
\item $t_{n0}\geq0$ for $\left\Vert T\right\Vert $-almost all points $p\in Supp(T)$,
that is $<\overrightarrow{T}(q),\mathbf{e}(q)>\geq0$, $\left\Vert T\right\Vert $-almost
everywhere; where $\mathbf{e}(q)$ is the unique horizontal (that
is, perpendicular to the fibers) $n$-plane at $q$ whose orientation
is preserved under $\pi_{*}$.
\item $\pi_{\#}(T)=1[M]$ as an $n$-dimensional current on $M$.
\item $\partial T=0$ (equivalently, for any $\Omega\subset M$, $\partial\left(T\rest\pi^{-1}(\Omega)\right)$
has support contained in $\partial\pi^{-1}(\Omega)$).
\end{enumerate}
\end{defn}
Note that each of these conditions is closed under weak convergence.
For the first, $t_{n,0}\geq0$ if $T(\phi)\geq0$ for all $\phi=\eta dx^{1}\wedge\cdots\wedge dx^{n}$,
where $\eta$ is a smooth, positive function with support in a neighborhood
of $p$. If $T_{i}$ is a sequence of such currents and $T_{i}\rightharpoonup T$,
then $T$ will also satisfy that condition. Similarly, for the second
condition, $\pi_{\#}(T)=1[M]$ if and only if $T(\pi^{*}(dV))=Vol(M)$,
which is again clearly closed under weak convergence. The third condition,
likewise, translates as $0=T(d\phi)$ for all smooth forms $\phi$,
which is also closed under weak convergence.

\begin{defn}
There is an $A>0$ so that the fiber bundle $B$ is contained in the
disk bundle $E_{A}\subset E$ defined by $E_{A}:=\left\{ v\in E\left|\left\Vert v\right\Vert <A\right.\right\} $,
by compactness of $B$. Define the space $\widetilde{\Gamma}(E)$
to be the set of all countably rectifiable, integer multiplicity,
$n$-dimensional currents which are quasi-sections in $E$, with support
contained in $E_{A }$, called \emph{(bounded)} \emph{rectifiable
sections} of $E$, which by the above is a weakly-closed set. The
space $\Gamma(E)$ of \emph{(strongly) rectifiable sections of} $E$
is the smallest sequentially weakly-closed space containing the graphs
of $C^{1}$ sections of $E$ which are supported within $E_{A}$.
\end{defn}
Thus, a quasi-section which is rectifiable and of integer multiplicity
is an element of $\widetilde{\Gamma}(E)$. It would seem to be a strictly
stronger condition for it to be in $\Gamma(E)$, however, it is shown
in \cite{Coventry} that, over a bounded domain $\Omega$, $\widetilde{\Gamma}(\Omega\times\R^{k})=\Gamma(\Omega\times\R^{k})$.
This extends to the statement that $\widetilde{\Gamma}(E)=\Gamma(E)$
for a vector bundle over a compact manifold $M$, since any such can
be decomposed into finitely many bounded domains where the bundle
structure is trivial, by a partition of unity argument.

The space $\widetilde{\Gamma}(B)$ of \emph{rectifiable sections}
of $B$ is the subset of $\widetilde{\Gamma}(E)$ of currents with
support in $B$, which is a weakly closed condition with respect to
weak convergence. Weak closure follows since, for any point $z$ outside
of $B$, there is a smooth form supported in a compact neighborhood
of $z$ disjoint from $B$. The space $\Gamma(B)$ of \emph{strongly}
rectifiable sections is the smallest sequentially, weakly-closed space
containing the graphs of $C^{1}$ sections of $B$. Since the fibers
of $B$ are compact, as is the base manifold $M$, minimal-mass elements
will exist in $\widetilde{\Gamma}(B)$ or $\Gamma(B)$, and mass-minimizing
sequences within any homology class will have convergent subsequences
in $\widetilde{\Gamma}(B)$ or $\Gamma(B)$. This follows from lower
semi-continuity with respect to convergence of currents, convexity
of the mass functional, and the closure and compactness theorems for
rectifiable currents. Closure of the conditions of definition (\ref{def:quasi-section})
under weak convergence will imply that the limits given by the closure
and compactness theorems, which are \emph{a priori} rectifiable currents,
are indeed rectifiable sections. For compact manifolds, as above,
$\widetilde{\Gamma}(E)=\Gamma(E),$ but it is \emph{not} the case
that $\widetilde{\Gamma}(B)=\Gamma(B)$ in general (see Proposition
(\ref{pro:Cart-not-cart}) below). 

\begin{rem}
A simple modification of the Federer-Flemming closure and compactness
theorems shows the following result: \cite{js-partial-regularity,js-1}
\end{rem}
\begin{prop}
Let $\{ T_{j}\}\subset\Gamma(B)$ (resp, $\widetilde{\Gamma}(B)$)
be a sequence with equibounded flat norm. Then, there is a subsequence
which converges weakly to a current $T$ in $\Gamma(B)$ (resp, $\widetilde{\Gamma}(B)$)
. 
\end{prop}
\begin{defn}
Given a current $T$, the induced measures $\left\Vert T\right\Vert $
and $\left\Vert T_{\alpha\beta}\right\Vert $ are defined locally
by:\begin{eqnarray*}
\left\Vert T_{\alpha\beta}\right\Vert (A) & := & \sup\left(T_{\alpha\beta}(\omega)\right),\,\textrm{and}\\
\left\Vert T\right\Vert (A) & := & \sup\left(\sum_{\alpha\beta}T_{\alpha\beta}(\omega)\right),\end{eqnarray*}
where the supremum in either case is taken over all $n$-forms on
$B$, $\omega\in E_{0}^{n}(B)$, with $comass(\omega)\leq1$ \cite[ 4.1.7]{GMT}
and $Supp(\omega)\subset A$.
\end{defn}

\subsection{\label{sub:Crofton's-formula}Crofton's formula}

The usual Crofton's formula (cf. for example \cite[3.2.26]{GMT})
for the measure of a rectifiable set states that, if $W$ is a rectifiable,
Hausdorff $n$-dimensional set in $\R^{n+k}$, then \[
\mathcal{H}^{n}(W)=\frac{1}{\beta(n+k,n)}\int_{p\in O^{*}(n+k,n)}\int_{\R^{n}}N(p|W,y)d\mathcal{L}^{n}(y)dV_{O*(N,n)}(p),\]
where $N(p|W,y)$ is the multiplicity at $y\in\R^{n}$ of the orthogonal
projection $p:\R^{n+k}\rightarrow\R^{n}$ restricted to $W$, $O^{*}(n+k,n)$
is the space of all such projections with the natural metric of total
volume 1, and $\beta(n+k,n)=\int_{p\in O^{*}(n+k,n)}\left\Vert p_{*}(P)\right\Vert dV_{O*(n+k,n)}(p)$. 

Since the mass of an integer-multiplicity, countably-rectifiable $n$-current
$T$ in $\R^{n+k}$ is the integral with respect to Hausdorff $n$-dimensional
measure restricted to the support of $T$ of the absolute value of
the multiplicity $\theta$, the mass of such a $T$ can be represented
by essentially the same integral-geometric formula. 

\begin{prop}
If $T$ is an integer-multiplicity, countably-rectifiable $n$-current
in $\R^{n+k}$, with multiplicity $\theta$, then the mass of $T$
is given by\[
\mathcal{M}(T)=\frac{1}{\beta(N,n)}\int_{p\in O^{*}(n+k,n)}\int_{\R^{n}}N(p|T,y,\theta)d\mathcal{L}^{n}(y)dV_{O*(N,n)}(p),\]
where $N(p|T,y,\theta)=N(p|Supp(T),y)|\theta|$ is the multiplicity
at $y\in\R^{n}$ of the orthogonal projection $p:\R^{N}\rightarrow\R^{n}$
restricted to $Supp(T)$, multiplied at each $z\in p^{-1}(y)\cap Supp(T)$
by $|\theta(z)|$, and $O^{*}(N,n)$ is the space of all such projections
with the natural metric of total volume 1. 
\end{prop}
For $T\in\widetilde{\Gamma}(B)$, and $i\in0,\ldots,n$, set $T_{i}=\sum_{|\beta|=i}T_{\alpha,\beta}$.
$T_{i}$ is the sum of the components of $T$ that have $i$ vertical
directions. Take $x_{0}\in M$ and $R>0$. Set $O^{*}(E,n,i)$ to
be the set of orthogonal projections from $\pi^{-1}(B(x_{0},R))\cong B(x_{0},R)\times\R^{k}\subset\R^{n+k}:=E$
which preserve $i$ vertical directions, that is, for which the kernel
contains an $\R^{k-i}$ inside of the fiber directions. Any such projection
is of course a direct product of projections $p_{1}:B(x_{0},R)\rightarrow\R^{n-i}$
and $p_{2}:\R^{k}\rightarrow\R^{i}$, so \[
O^{*}(E,n,i)=O^{*}(\R^{n},n-i)\times O^{*}(\R^{k},i).\]
 If $T$ is a smooth graph, $T=graph(u)$, then\[
\mathcal{M}(T_{i}\rest\pi^{-1}(B(x_{0},R)))=\int_{B(x_{0},R)}\left\Vert \nabla u^{\wedge i}\right\Vert dV.\]

\begin{prop}
\[
\mathcal{M}(T_{i}\rest\pi^{-1}(B(x_{0},R)))=\frac{\int_{p\in O^{*}(E,n,i)}\int_{\R^{n}}N(p|T,y,\theta)d\mathcal{L}^{n}(y)dV_{O*(E,n,i)}(p)}{\beta(n,n-i)\beta(k,i)}.\]

\end{prop}
\begin{proof}
(Compare \cite[3.16]{Morgan}) \begin{eqnarray*}
\mathcal{M}(T_{i}\rest\pi^{-1}(B(x_{0},R))) & = & \sup\left\{ \left.T_{i}(\phi)\right|comass(\phi)=1\right\} \\
 & = & \sup\left\{ \left.T(\phi)\right|comass(\phi)=1,\phi=\sum_{|\beta|=i}\phi_{\alpha\beta}dx^{\alpha}\wedge dy^{\beta}\right\} \\
 & = & \sup\left\{ \left.\int_{Supp(T)}<\overrightarrow{T}(z),\phi>\theta(z)d\left\Vert T\right\Vert \right|\begin{array}{l}
comass(\phi)=1,\\
\phi=\sum_{|\beta|=i}\phi_{\alpha\beta}dx^{\alpha}\wedge dy^{\beta}\end{array}\right\} \\
 & = & \frac{\int_{Supp(T)}\int_{p\in O^{*}(E,n,i)}\left\Vert p_{*}(\overrightarrow{T}(z))\right\Vert dV_{O*(E,n,i)}(p)\theta(z)d\left\Vert T\right\Vert }{\beta(n,n-i)\beta(k,i)}\\
 & = & \frac{\int_{p\in O^{*}(E,n,i)}\int_{Supp(T)}\left\Vert p_{*}(\overrightarrow{T}(z))\right\Vert \theta(z)d\left\Vert T\right\Vert dV_{O*(E,n,i)}(p)}{\beta(n,n-i)\beta(k,i)}\\
 & = & \frac{\int_{p\in O^{*}(E,n,i)}\int_{\R^{n}}N(p|T,y,\theta)d\mathcal{L}^{n}(y)dV_{O*(E,n,i)}(p)}{\beta(n,n-i)\beta(k,i)},\end{eqnarray*}
where the last step follows from the general area-coarea formula. 
\end{proof}

\section{Existence of horizontal cones\label{sec:Existence-of-horizontal}}

A current $C\in\widetilde{\Gamma}(B(x_{0},R)\times F)$ is an $h$-\emph{cone},
or a \emph{horizontal cone}, at $x_{0}$ if $(h_{\lambda})_{\#}(C)=C$.
From \cite{js-1}, a \emph{tangent} $h$-\emph{cone} at $x_{0}\in M$
of a rectifiable section $T\in\widetilde{\Gamma}(B)$ should be the
limit of horizontal dilations of $T$. First, restrict $T$ to $\pi^{-1}(B(x_{0},r_{0}))\approxeq B(x_{0},r_{0})\times F$.
Then, for $0<\lambda<r_{0}$, and $r>0$, set $h_{\lambda}:B(x_{0},\lambda r)\times F\rightarrow B(x_{0},r)\times F$
by $h_{\lambda}(x,v)=(x_{0}+(x-x_{0})/\lambda,v)$, and set $T_{\lambda}:=(h_{\lambda})_{\#}(T\rest B(x_{0},\lambda r)\times F)$.
In the case where $T=graph(u)$, then $T_{\lambda}$ is the graph
of $u_{\lambda}$ defined by $u_{\lambda}(x)=u(x_{0}+\lambda(x-x_{0}))$.
Then, for a sequence $\lambda_{i}\downarrow0$, the $h$-cone $H$
of $T$ at $x_{0}$ is the weak limit $H=\lim_{k}\left(h_{\lambda_{k}}\right)_{\#}(T\rest B(x_{0},\lambda_{k})\times F)$,
if that limit exists. Note that, as $\lambda\downarrow0$, the curvature
of the base will approach 0 and the bundle will become flat. The $h$-cone
is then defined on the Euclidean product $B(x_{0},r)\times F\subset\R^{n}\times F$. 

It was shown in \cite{js-1} that, for mass-minimizing rectifiable
sections as constructed in \cite{js-partial-regularity}, $h$-cones
always exist for some sequence of dilations, since a simple monotonicity
result shows that the set of dilations $T_{\lambda}$ will have equibounded
mass. We provide here a more direct proof of this fact in the case
we need. Note that the existence of $h$-cones is established only
for the mass-minimizing currents (with good partial-regularity) shown
to exist in \cite{js-partial-regularity}, which are limits of a sequence
of minimizers of functionals with an additional penalty term. It is
not known whether other mass-minimizers exist, without the required
partial regularity. 

For the moment, consider an arbitrary bundle $B\rightarrow M$ with
compact fiber $F$. Let $T$ be a {}``good'' mass-minimizing rectifiable
section, which is regular over an open dense subset. As before, set
$T_{i}=\sum_{|\beta|=i}T_{\alpha,\beta}$. From \cite[3.3.27]{GMT},
\[
\mathcal{M}(T\rest\pi^{-1}(B(x_{0},R))\leq\sum_{i=0}^{n}\mathcal{M}(T_{i}\rest\pi^{-1}(B(x_{0},R)).\]
This also follows directly from the triangle inequality.

In order to show that a sequence $(h_{\lambda})_{\#}(T\rest\pi^{-1}(B(x_{0},\lambda R)))$
of stretches converges, we need to show that each component $(h_{\lambda})_{\#}(T_{i}\rest\pi^{-1}(B(x_{0},\lambda R)))$
has mass bounded independently of $\lambda$. 

We use the result from \cite[Proposition 4.1]{js-partial-regularity},
stating that, since $T$ is mass-minimizing and is the limit of penalty-minimizers,
every point in $supp(T)$ has mass-density at least 1, and satisfies
standard monotonicity inequalities, $\mathcal{M}(T\rest B(z,\epsilon))\leq A\epsilon^{n}$. 

Consider $T\rest\pi^{-1}(B(x_{0},R))$. For each $z\in Supp(T\rest\pi^{-1}(B(x_{0},R)))$,
if $\epsilon>0$ is sufficiently small, the previous estimate holds
on $T\rest B(z,\epsilon)$, $\mathcal{M}(T\rest B(z,\epsilon))\leq A_{z}\epsilon^{n}$.
Since $Supp(T\rest\pi^{-1}(B(x_{0},\lambda R)))$ is compact, there
is a finite subcover $\mathcal{U}$ of such balls, with minimum radius
$\epsilon$. Let $A$ be the maximum of the constants $A_{z}$ for
these balls. Now, let $p\in O^{*}(E,n,i)$. Any ball centered at $y\in Im(p)\subset\R^{n}$
of radius $\epsilon$ will be such that $p^{-1}(y)$ meets finitely
many balls in this cover $\mathcal{U}$ (since the whole cover is
finite). The mass of the image of each of these balls is less than
the mass of the ball in $T$, since projection is mass-decreasing,
so the total image mass within that ball, counting multiplicities,
is less than the number of balls in the cover which intersect $p^{-1}(y)$,
times $A\epsilon^{n}$. Thus, there is a constant $C$ so that\[
\mathcal{M}(p_{\#}(T\rest\pi^{-1}(B(x_{0},R))))\leq C\mathcal{L}^{n}(p(Supp(T\rest\pi^{-1}(B(x_{0}R))))),\]
where $p_{\#}(T\rest\pi^{-1}(B(x_{0},R)))$ is the Crofton push-forward
current as in \S \ref{sub:Crofton's-formula}, with multiplicity
function $N(p|T,y,\theta)$ at each point in the image. 

Similarly, \begin{eqnarray*}
\mathcal{M}(T_{i}\rest\pi^{-1}(B(x_{0},R))) & \leq & C\mathcal{L}^{i}(p(F))\omega_{n-i}R^{n-i}\end{eqnarray*}
where $p(F)$ is the image of the fiber $F$ in $\R^{i}$ ($F$ is
a submanifold of $E_{x}\cong\R^{k}$), maximized over all $p\in O^{*}(E,n,i)$.
This inequality follows since the image of the projection of $T$
is contained in the image of $F\times B(x_{0},R)$. 

For precisely the same reasons, with the same constants, \begin{eqnarray*}
\mathcal{M}((h_{\lambda})_{\#}(T_{i}\rest\pi^{-1}(B(x_{0},\lambda R)))) & \leq & C\mathcal{L}^{i}(p(F))\omega_{n-i}R^{n-i},\end{eqnarray*}
since the factor of $\lambda$ coming from the stretch simply expands
the image of each projection until it again is contained within the
image of $F\times B(x_{0},R)$. The conclusion of this argument is
the following proposition:

\begin{prop}
\label{pro:limit_of_stretches}$\mathcal{M}((h_{\lambda})_{\#}(T\rest\pi^{-1}(B(x_{0},\lambda R))))$
is bounded, independently of $\lambda$. Thus, given a sequence $\lambda_{m}\downarrow0$,
a subsequence of $(h_{\lambda_{m}})_{\#}(T\rest\pi^{-1}(B(x_{0},\lambda_{m}R)))$
converges to a rectifiable section $T_{0}$ in $\widetilde{\Gamma}(B(x_{0},R)\times F)$. 
\end{prop}
\begin{proof}
Set $T^{m,R}:=(h_{\lambda_{m}})_{\#}(T\rest\pi^{-1}(B(x_{0},\lambda_{m}R)))$.
Then, by taking a diagonal subsequence, for each $j\in\Z$ there is
a current $T^{j}\in\widetilde{\Gamma}(B(x_{0},j)\times F)$ so that
$(h_{\lambda_{m}})_{\#}(T\rest\pi^{-1}(B(x_{0},\lambda_{m}j)))\rightharpoonup T^{j}$
and $T^{j}\rest B(x_{0},l)\times F=T^{l}$, whenever $j>l$, so that
there is a current $T^{0}$ on $R^{n}\times F$ which restricts to
each of these $T^{j}$. 
\end{proof}
We now specialize to the case of an $S^{n-1}$-bundle over a compact
$n$-manifold $M$. 

\begin{prop}
\label{pro:existence of h-cones}Let $B\rightarrow M$ be an $(n-1)$-sphere
bundle over a compact $n$-manifold $M$. Let $T$ be a good mass-minimizing
rectifiable section as before. Assume that $x_{0}\in M$ is a pole
point of $T$ so that the Hausdorff dimension of the pole is $(n-1)$,
that is, that the projection map $\phi_{r}:S^{(n-1)}(r)\times S^{n-1}\rightarrow S^{n-1}$,
inducing a Crofton projection $(\phi_{r})_{\#}(T\rest S^{(n-1)}(r)\times S^{n-1})\in\mathbf{R}^{n-1}(S^{n-1})$,
has limit having positive $(n-1)$-dimensional mass $A$ for some
subsequence of the sequence $r_{m}=\lambda_{m}R$. Then the current
$T^{0}$ of Proposition (\ref{pro:limit_of_stretches}) will be an
$h$-cone. 
\end{prop}
\begin{proof}
Since each $T^{j}$ minimizes the scaled and stretched functional
\[
\mathcal{V}^{j}(S):=\mathcal{V}((h_{\lambda_{j}}^{-1})_{\#}(S))/\left(\lambda_{j}R\mathcal{M}((\phi_{\lambda_{j}R})_{\#}(T\rest S^{n-1}(\lambda_{j}R)\times S^{n-1})\right),\]
 $T^{0}$ will minimize the limiting functional \begin{eqnarray*}
\mathcal{V}_{0}(S\rest B^{n}(x_{0},R)\times S^{n-1}) & = & \lim\mathcal{V}^{j}(S\rest B^{n}(x_{0},R)\times S^{n-1})\\
 & = & \mathcal{M}(S_{n-1}\rest B^{n}(x_{0},R)\times S^{n-1}),\end{eqnarray*}
where $S_{n-1}:=\sum_{|\alpha|=1}S_{\alpha\beta}$ is that part of
the current $S$ which has $(n-1)$ vertical components, one horizontal
component. The stretched functionals $\mathcal{V}^{j}$, as $j\rightarrow\infty$,
magnify the terms with more vertical components by the effect of $(h_{\lambda_{j}}^{-1})_{\#}$,
and under the assumption that the pole at $x_{0}$ has Hausdorff dimension
$(n-1)$ that highest-order term will dominate all others in the normalized
limit. This reduces to \[
\int_{B(x_{0},R)}\left\Vert \nabla u^{\wedge(n-1)}\right\Vert dV,\]
if $S$ is a smooth graph $S=graph(u)$. Note also that $\mathcal{M}^{n-1}(S_{n-1}\rest\partial B^{n}(x_{0},R)\times S^{n-1})$
is the $(n-1)$-dimensional mass of the projection $(\phi_{R})_{\#}(T^{0}\rest S^{n-1}(R)\times S^{n-1})$.
Since $T^{0}$ minimizes, for any $R$\begin{eqnarray*}
\mathcal{V}_{0}(T^{0}\rest B^{n}(x_{0},R)\times S^{n-1}) & \leq & \mathcal{V}_{0}(C(T^{0}\rest\partial B^{n}(x_{0},R)\times S^{n-1}))\\
 & = & R\mathcal{M}^{n-1}(T_{n-1}^{0}\rest\partial B^{n}(x_{0},R)\times S^{n-1})\\
 & = & R\mathcal{M}^{n-1}\left((\phi_{R})_{\#}(T^{0}\rest S^{n-1}(R)\times S^{n-1})\right),\end{eqnarray*}
where $C()$ denotes the $h$-cone over the boundary $T_{0}\rest\partial B^{n}(x_{0},R)\times S^{n-1})$.
On the other hand, by slicing\begin{eqnarray*}
\frac{d}{dR}\mathcal{V}_{0}(T^{0}\rest B^{n}(x_{0},R)\times S^{n-1}) & \geq & \mathcal{M}^{n-1}(T_{n-1}^{0}\rest\partial B^{n}(x_{0},R)\times S^{n-1})\\
 & = & \mathcal{M}^{n-1}\left((\phi_{R})_{\#}(T^{0}\rest S^{n-1}(R)\times S^{n-1})\right),\end{eqnarray*}
so that \begin{eqnarray*}
 &  & \frac{d}{dR}\frac{\mathcal{V}_{0}(T^{0}\rest B^{n}(x_{0},R)\times S^{n-1})}{R}\\
 & = & \frac{\frac{d}{dR}\left(\mathcal{V}_{0}(T^{0}\rest B^{n}(x_{0},R)\times S^{n-1})\right)R-\mathcal{V}_{0}(T^{0}\rest B^{n}(x_{0},R)\times S^{n-1})}{R}\\
 & \geq & \frac{\mathcal{M}^{n-1}\left((\phi_{R})_{\#}(T^{0}\rest S^{n-1}(R)\times S^{n-1})\right)R-\mathcal{V}_{0}(T^{0}\rest B^{n}(x_{0},R)\times S^{n-1})}{R}\\
 & \geq & 0,\end{eqnarray*}
and so $\mathcal{V}_{0}(T^{0}\rest B^{n}(x_{0},R)\times S^{n-1})/R$
is an increasing function of $R$. However, since $T^{0}$ is invariant
at least under the sequence of stretches by $h_{\lambda_{j}}$, the
projected mass $\mathcal{M}^{n-1}\left((\phi_{R})_{\#}(T^{0}\rest S^{n-1}(R)\times S^{n-1})\right)$
must be the same for $R=\lambda_{j}R_{0}$, so that the values of
$\mathcal{M}^{n-1}\left((\phi_{R})_{\#}(T^{0}\rest S^{n-1}(R)\times S^{n-1})\right)$
repeat over the intervals $[\lambda_{j+1}R_{0},\lambda_{j}R_{0}]$
and the increasing function $\mathcal{V}_{0}(T^{0}\rest B^{n}(x_{0},R)\times S^{n-1})/R$
satisfies\[
\frac{\mathcal{V}_{0}(T^{0}\rest B^{n}(x_{0},R)\times S^{n-1})}{R}\leq\mathcal{M}^{n-1}\left((\phi_{R})_{\#}(T^{0}\rest S^{n-1}(R)\times S^{n-1})\right)\]
so \[
\mathcal{V}_{0}(T^{0}\rest B^{n}(x_{0},R)\times S^{n-1})/R\leq\inf\left(\mathcal{M}^{n-1}\left((\phi_{R})_{\#}(T^{0}\rest S^{n-1}(R)\times S^{n-1})\right)\right).\]
 However, since \begin{eqnarray*}
\mathcal{V}_{0}(T^{0}\rest B^{n}(x_{0},R)\times S^{n-1}) & \ge & \int_{0}^{R}\mathcal{M}^{n-1}\left((\phi_{r})_{\#}(T^{0}\rest S^{n-1}(r)\times S^{n-1})\right)dr\\
 & \geq & R\inf\left(\mathcal{M}^{n-1}\left((\phi_{R})_{\#}(T^{0}\rest S^{n-1}(R)\times S^{n-1})\right)\right),\end{eqnarray*}
all of these inequalities must be equalities, and necessarily $\mathcal{M}^{n-1}\left((\phi_{R})_{\#}(T^{0}\rest S^{n-1}(R)\times S^{n-1})\right)$
must be constant. Moreover, \begin{eqnarray*}
\mathcal{V}_{0}(T^{0}\rest B^{n}(x_{0},R)\times S^{n-1}) & = & \mathcal{V}_{0}(C(T^{0}\rest\partial B^{n}(x_{0},R)\times S^{n-1})),\end{eqnarray*}
and since any change with respect to the radial direction (of positive
measure) would introduce a strict inequality in that integral, $T^{0}\rest B^{n}(x_{0},R)\times S^{n-1}=C(T^{0}\rest\partial B^{n}(x_{0},R)\times S^{n-1})$
$T^{0}$-almost everywhere. Thus $T^{0}$ is an $h$-cone. 
\end{proof}
\begin{thm}
\label{thm:Cart=3Dcart-if-poles}If $B$ is an $(n-1)$-sphere bundle
over a compact $n$-manifold $M$, and if $T\in\widetilde{\Gamma}(B)$
is a smooth graph except on a finite set of fibers $\pi^{-1}(x_{i})$,
so that the degree of each singular fiber is $0$ and so that there
is an $h$-cone at each fiber, then $T\in\Gamma(B)$.
\end{thm}
\begin{proof}
The only part of this statement requiring proof is that, in a neighborhood
of each singular fiber, the current is a limit of smooth currents.
Certainly, if the degree of any of the singular fibers is nonzero
it cannot be in $\Gamma(B)$. If the degree is 0, however, since the
graph is smooth within the boundary spheres $S^{n-1}(r)\times S^{n-1}$
of $B\rest\pi^{-1}(B(x_{0},r))\cong B(x_{0},r)\times S^{n-1}$, the
$h$-cone is a cone over a current $S\in\Gamma(S^{n-1}(1)\times S^{n-1})$,
in fact, $S$ is the limit of the smooth sequence of stretches of
$T\rest S^{n-1}(r)\times S^{n-1}$. Since the degree of the singularity
is 0, each graph $T\rest S^{n-1}(r)\times S^{n-1}$ is (smoothly)
homotopic to the constant map, mapping $S^{n-1}(r)$ to $p_{0}\in S^{n-1}$.
If $H(x,t):S^{n-1}(r)\times[0,1]\rightarrow S^{n-1}$ is that homotopy,
then the graph $G(y):B(x_{0},r)\rightarrow S^{n-1}$ defined by $G(y)=H(ry/|y|,1-|y|/r)$
will be a smoothable graph which can be extended to a section of $B$
agreeing with $T$ outside of this neighborhood. Clearly, given a
sequence $r_{i}\rightarrow0$, the maps \[
T_{i}=\begin{cases}
G_{r_{i}}, & d(x_{0},x)<r_{i}\\
T, & d(x_{0},x)\geq r_{i}\end{cases}\]
will be a sequence of currents converging weakly to $T$, which are
smooth in a neighborhood of the pole point $x_{0}$. Since there are
finitely many singular points of $T$ by hypothesis, iterating this
construction will generate a sequence of smooth currents converging
to $T$. 
\end{proof}

\subsection{Rectifiable foliations, rectifiable sections.}

Consider now the case where $B$ is the subbundle $T_{1}(M)$ of unit
vectors in $T_{*}(M)$. The connection used to define the metric on
$T_{*}(M)$ restricts to an associated connection on $T_{1}(M)$,
since the connection is a metric connection, and defines a metric
on $T_{1}(M)$ as before. 

\textit{Rectifiable 1-dimensional foliations} on $M$ are rectifiable
sections of $T_{*}(M)$ whose support lies within $T_{1}(M)$. As
above, this condition will be weakly closed, so that the Federer closure
and compactness theorems hold.

\begin{thm}
\label{basic regularity} \cite{js-1,js-partial-regularity} For any
homology class of sections in $\widetilde{\Gamma}(T_{1}(M))$, there
is a mass-minimizing rectifiable foliation $\mathcal{F}$ with support
which is the Gauss map of a $C^{1}$ graph over an open, dense subset
of $M$. 
\end{thm}
The \textit{regular points} of a rectifiable foliation $S$ correspond
to points where the Gauss map is continuous, and singularities, or
\textit{pole points}, are points $x\in M$ where the Gauss map is
discontinuous. Equivalently, pole points are those $x\in M$ for which
the set $Supp(S)\cap\pi^{-1}(x)$ consists of more than one point.
Points of $Supp(S)$ lying over pole points are called \textit{pole
elements}.

\section{\label{sec:-h--Cones}The degree of a pole point}

Let $S$ be a mass-minimizing rectifiable section of $T_{1}(M^{3})$.
By \cite{dim-3}, there is such a minimizer with only a finite number
of pole points, each of which contains the entire fiber in the support
of $S$. Note that, in\cite{dim-3}, the results need to be modified
to indicate that \cite{js-partial-regularity} does not show that
any minimizer has the required smoothness, only that there is one
minimizer with the claimed partial regularity. Assume that $S$ is
such a minimizer. The question of regularity of a mass-minimizing
section of $T_{1}(M^{3})$ becomes whether such a pole point can exist. 

If $x_{0}$ is an isolated pole point of $S\in\widetilde{\Gamma}(T_{1}(M^{n}))$,
by Proposition (\ref{pro:existence of h-cones}) there is an $h$-cone
centered at $x_{0}$. By Theorem (\ref{thm:Cart=3Dcart-if-poles}),
$S\in\Gamma(T_{1}(M^{n}))$. In addition, the $h$-cone at $x_{0}$
is a rectifiable section $\psi\in\Gamma(B^{n}\times S^{n-1})$, when
restricted to a ball of radius $1$ in the base. Slicing the $h$-cone
$\psi$ by a cylinder of radius $r$ generates a rectifiable current
$C$ in $S^{n-1}(r)\times S^{n-1}\cong S^{n-1}\times S^{n-1}$ for
almost any $r$ by slicing theory. Since $\psi$ is an $h$-cone,
however, $C$ is independent of $r$, thus the slice is rectifiable
for all $r$, and so is in $\Gamma(S^{n-1}\times S^{n-1})$ as a bundle
over the first factor. The key to existence of such a singularity
is the degree of the current $C$.

Since $S$ has no interior boundary, and by \cite{dim-3} the support
of $S$ contains all of $\pi^{-1}(x_{0})$, the \textit{image} $I(C)$
of $C$, defined as the push-forward image $(\Pi_{2})_{\#}(C)$, for
$\Pi_{2}:S^{n-1}\times S^{n-1}\rightarrow S^{n-1}$ the projection
onto the second factor (the fiber), must have support the entire sphere.
The degree of a rectifiable section $C\in\widetilde{\Gamma}(S^{n-1}\times S^{n-1})$
is defined by \[
deg(C):=\int_{I(C)}dV_{S^{n-1}}=\int_{C}\Pi_{2}^{*}(dV_{S^{n-1}})=C(\Pi_{2}^{*}(dV_{S^{n-1}})),\]
which is clearly a weakly closed condition. If $C$ is the graph of
a smooth map $\phi:S^{n-1}\rightarrow S^{n-1}$, then $deg(C)=deg(\phi)$,
and in particular, if $\phi$ is the restriction of a smooth map $\Phi:B^{n}\rightarrow S^{n-1}$
to the boundary, then $deg(C)=0$. By taking transfinite limits, if
$C$ arises from the $h$-cone of a strongly rectifiable section $S\in\Gamma(B^{n}\times S^{n-1})$,
$deg(C)=0$ since $C$ is a weak limit of degree-zero currents.

\begin{defn}
The \emph{degree} of a pole point $x_{0}\in M^{n}$ of a rectifiable
section $S\in\Gamma(B)$, where $B$ is an $S^{n-1}$-bundle over
$M$, $deg(S,x_{0})$, is the degree of the restriction of an $h$-cone
$\psi$ of $S$ to the boundary $\psi\rest S^{n-1}(r)\times S^{n-1}$.
As mentioned earlier, if $S\in\Gamma(B)$, $deg(S,x_{0})=0$. 

We now return to the claim in Section \ref{sec:Rectifiable-Sections.}
that not all weak rectifiable sections are strong rectifiable sections.
\end{defn}
\begin{prop}
\label{pro:Cart-not-cart}$\widetilde{\Gamma}(T_{1}(S^{2}))\neq\Gamma(T_{1}(S^{2}))$.
\end{prop}
\begin{proof}
Since there are no continuous sections of $T_{1}(S^{2})$, that is,
$\Gamma(T_{1}(S^{2}))=\emptyset$, it suffices to show that $\widetilde{\Gamma}(T_{1}(S^{2}))\neq\emptyset$.
Given a point $p\in S^{2}$, and $v\in T_{1}(S^{2},p)$, translate
$v$ parallel to itself along longitudes to $-p$. The rectifiable
section generated by this procedure will have a singular point at
$-p$, with the entire fiber of the sphere bundle in the support over
$-p$. Since there is no boundary and it projects to $1[S^{2}]$ on
$S^{2}\backslash\{-p\}$, it is an element of $\widetilde{\Gamma}(T_{1}(S^{2}))$. 
\end{proof}
\begin{rem}
Of course, this current is an element of $\Gamma(T_{*}(S^{2}))$,
and is the limit of a sequence of smooth vector fields, each of which
has a zero of degree 2 at $-p$, with length 1 outside of neighborhoods
of $-p$. It should also be noted that this topological obstruction
is not the only way that it can be possible for $\widetilde{\Gamma}(B)\neq\Gamma(B)$
for $B$ an $(n-1)$-sphere bundle on an $n$-manifold. Since the
degree of an isolated singularity is local, it follows that any isolated
singularities of $T\in\Gamma(B)$ will have degree 0. But even on
a sphere bundle $B$ with global smooth sections, it is easy to construct
singular sections with two isolated singularities, one of degree 2
and the other of degree -2. Such singular sections are clearly in
$\widetilde{\Gamma}(B)\backslash\Gamma(B)$.
\end{rem}

\section{Non-Existence of Isolated Singularities}

Now that we have shown that an isolated pole of a volume-minimizing
section $S$ of $T_{1}(M)$ necessarily stretches to an energy-minimizing
section $S_{0}$ for the limiting volume $\mathcal{V}_{0}$, we proceed
to show that it cannot exist if the degree of the pole is $0$. As
before, set $S_{0}\in\widetilde{\Gamma}_{R}(B(0,r)\times S^{n-1})$
to be an $h$-cone of $S$, and set $C:=S_{0}\rest S^{n-1}\times S^{n-1}$
($r$ may be assumed to be larger than 1).

\begin{thm}
\label{no-isolated}There is a volume minimizing section $S$ of $T_{1}(M)$
with no degree-zero isolated pole points $x_{0}$, with $supp(S)\cap\pi^{-1}(x_{0})=S^{n-1}=\pi^{-1}(x_{0})$.
\end{thm}
\begin{proof}
Let $S$ be a mass-minimizing section which is continuous over an
open, dense subset, as guaranteed by \cite{js-partial-regularity},
as discussed in the previous section. Assume that $S$ has a degree-zero
isolated singularity $x_{0}$, with the entire fiber contained in
the support of $S$. There is a $h$-cone $S_{0}$ of $S$ at $x_{0}$
by \S\ref{sec:Existence-of-horizontal}. the current $C=S_{0}\rest S^{n-1}(r)\times S^{n-1}$
has degree 0, as in \S\ref{sec:-h--Cones}, and so there is a rectifiable
current $F$ so that $\partial F=C-graph(constant)$ in $S^{n-1}\times S^{n-1}$.
In fact, the $h$-cone $S_{0}$ can be used to construct such a current
$F_{0}$ which is a {}``rectifiable homotopy'', that is, which extends
to a rectifiable section on $\left(S^{n-1}\times I\right)\times S^{n-1}$as
an $n$-dimensional current with $\partial F=C\times0-S^{2n-1}\times\{ pt\}\times1$.
For each $i$ in a sequence $S_{i}\in\widetilde{\Gamma}(B(x_{0},1)\times S^{n-1})$
converging to the $h$-cone $S_{0}$, and for each $r$, $S_{i}\rest\partial B(x_{0},r)\times S^{n-1}=S_{i}(r)\rest S^{n-1}\times S^{
}$ is a smooth graph of degree 0, so there is a rectifiable current
{}``fence'' $F_{i}(r)$ of dimension $n$ so that $\partial F_{i}(r)=S_{i}(r)\rest S^{n-1}\times S^{
}-graph(constant)\rest S^{n-1}\times S^{
}$. Since $S_{i}\rightharpoonup S_{0}$, which is a cone, $F_{i}(r)$
can be chosen with bounded mass, so there is a convergent subsequence
with limit $F_{0}(r)$. Since $S_{0}$ is an $h$-cone, it may be
assumed that $F_{0}(r)=(h_{r})_{\#}(F_{0}(1))$. The current \[
S_{r}:=S_{0}\rest B(x_{0},R)\backslash B(x_{0},r)\times S^{n-1}+F_{0}(r)+graph(constant)\rest B(x_{0},r)\]
 has the same boundary as $S_{0}\rest B(x_{0},R)$. However, $S_{0}\rest B(x_{0},R)$
minimizes the limiting functional $\mathcal{V}_{0}$, so, independent
of $r$, $\mathcal{V}_{0}(S_{r})\geq\mathcal{V}_{0}(S_{0})$. But,
\begin{eqnarray*}
\mathcal{V}_{0}(S_{r}) & = & \mathcal{V}_{0}\left(S_{0}\rest B(x_{0},R)\backslash B(x_{0},r)\times S^{n-1}+F_{0}(r)+graph(constant)\rest B(x_{0},r)\right)\\
 & = & \mathcal{V}_{0}\left(S_{0}\rest B(x_{0},R)\backslash B(x_{0},r)\times S^{n-1}\right)+\mathcal{V}_{0}\left(F_{0}(r)\right),\end{eqnarray*}
since $\mathcal{V}_{0}(graph(constant))=0$. In addition, $\mathcal{V}_{0}\left(F_{0}(r)\right)=Ar^{n-1}$
since only the $n-1$ base dimensions are stretching with $r$. However,
\[
\mathcal{V}_{0}(S_{0}\rest B(x_{0},R))-\mathcal{V}_{0}\left(S_{0}\rest B(x_{0},R)\backslash B(x_{0},r)\times S^{n-1}\right)=Br^{n}\]
The constants $A$ and $B$ do not depend upon $R$, except for the
limitation that $r<R$. Clearly, for $R$ sufficiently large $\mathcal{V}_{0}(S_{0}\rest B(x_{0},R))-\mathcal{V}_{0}(S_{r})=Br^{n}-Ar^{n-1}$
will eventually be positive for some $r$ large enough, contradicting
the fact that $S_{0}\rest B(x_{0},R)$ minimizes $\mathcal{V}_{0}$
there.
\end{proof}
\begin{cor}
If $M$ is a compact $3$-manifold, then there is a volume-minimizing
one-dimensional foliation of class $C^{1}$.
\end{cor}
\begin{proof}
By \cite{dim-3}, there is a volume-minimizing rectifiable section
of $T_{1}(M)$ with only isolated singular points, for which the support
of each contains the entire fiber. Such isolated poles cannot exist
by the theorem, so there is a rectifiable section with no poles, so
that the section is continuous on all of $M$. Since that section
is the tangent field of the foliation, the foliation is of class $C^{1}$.
\end{proof}
\begin{defn}
Sharon Pedersen, in \cite{Pedersen}, defined, for each $n\geq1$,
a rectifiable section $P_{n}$ of $T_{1}(S^{2n+1})$, defined by parallel
translation of a unit vector $v\in T_{1}(S^{2n+1},x)$ along meridians
to $-x$. This is a rectifiable foliation and a minimal submanifold
except over a single point, and was shown to have, for $n>1$, much
smaller volume than the foliations defined by the standard Hopf fibrations.
She conjectured that this current might minimize volume amongst rectifiable
sections of $T_{1}(S^{2n+1})$, but this is not the case as shown
below.
\end{defn}
\begin{cor}
The rectifiable sections $P_{n}$ of $T_{1}(S^{2n+1})$ are not volume-minimizing
rectifiable foliations.
\end{cor}
\begin{proof}
The singularity at $-x$ of such a foliation is precisely the kind
shown to not exist by Theorem (\ref{no-isolated}).
\end{proof}


\begin{thebibliography}{10}
\bibitem{Almgren}F. J. Almgren, Jr. \textit{Q-valued functions minimizing Dirichlet's
integral and the regularity of area minimizing rectifiable currents
up to codimension two ,} Bull. (New Series) Am. Math. Soc., \textbf{8}
(1983), 327--328. 
\bibitem{Coventry}A. Coventry, The graphs of smooth functions are dense in the space
of Cartesian currents on a smooth bounded domain, but not every Cartesian
current is the limit of such graphs, Research reports (Mathematics),
number CMA-MRR 45- 98, Australian National University Publications,
1998.
\bibitem{GMT}Herbert Federer, \textit{Geometric Measure Theory,} Springer-Verlag
1969. 
\bibitem{GZ}Herman Gluck and Wolfgang Ziller, \textit{On the volume of a unit
vector field on the three-sphere,} Commentarii Mathematici Helvitici,
\textbf{61} (1986), 177--192. 
\bibitem{j}D. L. Johnson, \textit{K\"{a}hler submersions and holomorphic connections,}
Jour. Diff. Geo. \textbf{15} (1980) , 71--79.
\bibitem{js-partial-regularity}D. L. Johnson and P. Smith\textit{, Partial regularity of mass-minimizing
rectifiable sections,} submitted for publication, http://arXiv.org/abs/math/0403483. 
\bibitem{js-1}D. L. Johnson and P. Smith, \textit{Regularity of volume-minimizing
graphs,} Indiana University Mathematics Journal, \textbf{44} (1995),
45--85. 
\bibitem{dim-3}D. L. Johnson and P. Smith, \textit{Regularity of mass-minimizing
one-dimensional foliations,} Analysis and Geometry on Foliated Manifolds,
Proceedings of the VII International Colloquium on Differential Geometry,
(1994), World Scientific, 81--98.
\bibitem{Lawlor}Gary Lawler,  \textit{A sufficient condition for a cone to be area-minimizing,}
Memoirs of the American Mathematical Society \textbf{91} (1991). 
\bibitem{MS}J. D. Moore and R. Schlafly, \textit{On equivariant isometric embeddings,}
Mathematische Zeitschrift \textbf{173} (1980), 119--133. 
\bibitem{Morgan}F. Morgan, \emph{Geometric Measure Theory, A Beginner's Guide}, Academic
Press, 1988; second edition, 1995.
\bibitem{Pedersen}Sharon Pedersen \textit{Volumes of vector fields on spheres,} \textbf{336}
(1993) Trans. Am. Math. Soc., 69--78. 
\bibitem{Sasaki}Takeshi Sasaki, \textit{On the differential geometry of tangent bundles
of Riemannian manifolds} \textbf{10} (1958), T\^ ohoku Math. J.,
338--354.\end{thebibliography}
\end{document}